\documentclass[10pt]{amsart}
\usepackage{amscd}
\usepackage{amssymb}
\usepackage{latexsym}
\usepackage{graphicx}

\newcommand{\embeds}{\hookrightarrow} 

\newcommand{\init}{\vert_{t = 0}}

\newcommand{\abs}[1]{\left\vert #1 \right\vert}
\newcommand{\fixedabs}[1]{\vert #1 \vert}
\newcommand{\bigabs}[1]{\bigl\vert #1 \bigr\vert}

\newcommand{\norm}[1]{\left\Vert #1 \right\Vert}
\newcommand{\fixednorm}[1]{\Vert #1 \Vert}
\newcommand{\bignorm}[1]{\bigl\Vert #1 \bigr\Vert}

\newcommand{\C}{\mathbb{C}}
\newcommand{\R}{\mathbb{R}}

\newcommand{\innerprod}[2]{\left\langle \, #1 , #2 \,
\right\rangle}

\newcommand{\angles}[1]{\langle #1 \rangle}

\DeclareMathOperator{\sgn}{sgn}

\newtheorem{theorem}{Theorem}
\newtheorem{lemma}{Lemma}
\newtheorem{corollary}{Corollary}

\theoremstyle{definition}

\theoremstyle{remark}
\newtheorem{remark}{Remark}

\title[1D-DKG]{Low regularity well-posedness of the Dirac-Klein-Gordon equations in one space dimension}

\author{Sigmund Selberg}
\address{Department of Mathematical Sciences\\ Norwegian University of Science and Technology\\ Alfred Getz' vei 1\\ N-7491 Trondheim\\ Norway}
\email{sselberg@math.ntnu.no}
\urladdr{www.math.ntnu.no/~sselberg}
\thanks{Both authors supported by Research Council of Norway, project 160192/V30, PDE and Harmonic Analysis}

\author{Achenef Tesfahun}
\address{Department of Mathematical Sciences\\ Norwegian University of Science and Technology\\ Alfred Getz' vei 1\\ N-7491 Trondheim\\ Norway}
\email{tesfahun@math.ntnu.no}

\subjclass[2000]{35Q40; 35L70}

\begin{document}

\begin{abstract} 
We extend recent results of S.\ Machihara and H.\ Pecher on
low regularity well-posedness of the Dirac-Klein-Gordon (DKG) system in one dimension. Our proof, like that of Pecher, relies on the null structure of DKG, recently completed by D'Ancona, Foschi and Selberg, but we show that in 1d the argument can be simplified by modifying the choice of projections for the Dirac operator. We also show that the result is best possible up to endpoint cases, if one iterates in Bourgain-Klainerman-Machedon spaces.
\end{abstract}

\maketitle

\section{Introduction}\label{Section1}

We consider the Dirac-Klein-Gordon system (DKG) in one space dimension,
\begin{equation}\label{DKG1}
\left\{
\begin{alignedat}{2}
&D_t \psi + \alpha D_x  \psi + M\psi= \phi \beta\psi,&
\qquad\qquad
&\left(D_t=-i \partial_t, \,\, D_x=-i \partial_x\right)
\\
&- \square \phi + m^2 \phi= \innerprod{\beta\psi}{\psi}_{\C^2},
&
&\left(\square = -\partial_t^2 + \partial_x^2 \right)
\end{alignedat}
\right.
\end{equation}
with initial data
\begin{equation}\label{data}
  \psi \init = \psi_0 \in H^s, \qquad \phi \init =
\phi_0 \in H^r, \qquad
\partial_t \phi \init = \phi_1 \in H^{r-1},
\end{equation}
where $\phi(t,x)$ is real-valued and $\psi(t,x) \in \C^2$ is the Dirac spinor, regarded as a column vector with components $\psi_1$, $\psi_2$; $M,m \ge 0$ are constants. The $2 \times 2$ matrices $\alpha$, $\beta$ should be hermitian and satisfy $\beta^2 = \alpha^2 = I$, $\alpha \beta + \beta \alpha = 0$. A particular representation is
$$
  \alpha =    \begin{pmatrix}
    0 & 1  \\
    1 & 0
  \end{pmatrix},
  \qquad
  \beta =  \begin{pmatrix}
    1 & 0  \\
    0 & -1
  \end{pmatrix}.
$$

Global well-posedness for DKG in 1d was proved by Chadam \cite{Chadam:1973}, for data \eqref{data} with $(r,s)=(1,1)$. Several authors have improved Chadam's result, in the sense that the required regularity $(r,s)$ has been lowered; see Table \ref{Table0} for an overview.

\begin{table}
\caption{Global well-posedness for \eqref{DKG1}, \eqref{data}}
\label{Table0}
\def\arraystretch{1.3}
\begin{center}
\begin{tabular}{|r|c|c|}
\hline
& $s$ & $r$
\\
\hline
Chadam \cite{Chadam:1973}, 1973 & 1 & 1
\\
\hline
Bournaveas \cite{Bournaveas:2000}, 2000 & 0 & 1
\\
\hline
Fang \cite{Fang:2004b}, 2004 & 0 & (1/2,1]
\\
\hline
Bournaveas and Gibbeson \cite{BournaveasGibbeson:2006}, 2006 & 0 & [1/4,1]
\\
\hline
Machihara \cite{Machihara:2006}, Pecher \cite{Pecher:2006}, 2006 & 0 & (0,1]
\\
\hline
\end{tabular}
\end{center}
\end{table}

The global results are obtained by first proving local well-posedness and then using the conservation of the charge norm $\fixednorm{\psi(t)}_{L^2}$ together with a suitable a priori estimate for $\phi(t)$, to show that the solution extends globally.

Thus, the main step is to prove local well-posedness, and the best such results to date are due to  Machihara \cite{Machihara:2006} and Pecher \cite{Pecher:2006}, who worked independently of each other. Machihara proved local well posedness of \eqref{DKG1} for data \eqref{data} with $(s,r)$ in the region
$$
  -\frac{1}{4} < s \le 0,
  \qquad
  2\abs{s}
  \le r,
  \qquad
  r \le 1+2s.
$$
Pecher obtained the region
$$
  s > -\frac{1}{4}, \qquad r>0,
  \qquad
  \abs{s}
  \le r,
  \qquad
  r < 1+2s,
  \qquad r \le 1+s.
$$
To compare the two results, note that in Pecher's region, intersected with the strip
$-1/4 < s \le 0$, the lower bound for $r$ is $\abs{s}$, which is better than Machihara's
lower bound $\abs{2s}$, but on the other hand, Pecher has $r < 1+2s$ in this strip, whereas Machihara has $r \le 1+2s$.

Here we prove local well-posedness in a strictly larger region of the $(s,r)$-plane, which contains the union of the Pecher's and Machihara's regions. In fact, we show that in the strip $-1/4 < s \le 0$, the bound $r \le 1+2s$ can be relaxed to $r \le 1+s$.

\begin{theorem}\label{Thm1}
The DKG system \eqref{DKG1} is locally well posed for data
\eqref{data} with \ $(s,r)$ in the region
\begin{equation*}
  s > -\frac{1}{4}, \qquad r>0,
  \qquad
  \abs{s}
  \le r \le 1+s.
 \end{equation*}
\end{theorem}

Moreover, we show that this result is best possible, except possibly for the endpoint $(s,r) = (0,0)$, if one uses iteration in Bourgain-Klainerman-Machedon spaces; see Section \ref{Optimality}.

Our proof of Theorem \ref{Thm1}, like Pecher's original proof, relies on the null structure of DKG, which was completed recently by D'Ancona, Foschi and Selberg \cite{Selberg:2006b}. To see the null structure, one starts by decomposing the spinor into eigenvectors of the Dirac operator. This approach was used by Beals and Bezard \cite{Beals:1996} to show that $\innerprod{\beta \psi}{\psi}$ is a null form.\footnote{The fact that this expression is a null form was proved even earlier by Klainerman and Machedon \cite{Klainerman:1994a}, but they used a different, more indirect method.} The new idea introduced in \cite{Selberg:2006b} is that this null form then appears again in the Dirac equation, after a duality argument. The null structure was used in \cite{Selberg:2006b} to prove almost optimal local well-posedness of the 3d DKG system, and in \cite{Selberg:2006c} to treat the 2d case. Pecher's proof for the 1d case follows closely the argument in \cite{Selberg:2006b}, but here we show that in 1d the argument can be simplified by choosing the Dirac projections in a different way.

This paper is organized as follows: In the next section we reduce Theorem \ref{Thm1} to two bilinear estimates, and introduce the main tools needed for their proofs, which are given in Section \ref{Thm1Proof}. In Section \ref{Optimality} we prove the optimality of our result, by constructing explicit counterexamples for the iterative estimates. In Section \ref{Thm2Proof} we prove a product law for Wave-Sobolev spaces (see Theorem \ref{Thm2}) which is needed for the proof of Theorem \ref{Thm1}.

Let us fix some notation. We use $\lesssim$ to mean $\le$ up to multiplication by a positive constant $C$ which may depend on $s$ and $r$. If $a,b$ are nonnegative quantities, $a \sim b$ means $b \lesssim a \lesssim b$. The Fourier transforms in space and space-time are defined by
\begin{align*}
  \widehat f(\xi) &= \int_{\R} e^{-ix\xi}f(x) \, dx,
  \\
  \widetilde u(\tau,\xi) &= \int_{\R^{1+1}} e^{-i(t\tau+x\xi)}
u(t,x) \, dt \, dx,
\end{align*}
so $\widetilde{D_x u} = \xi \widetilde{u}$, $\widetilde{D_t
u} = \tau \widetilde{ u}$. $H^s = H^s(\R)$ is the Sobolev space with norm
$$
  \norm{f}_{H^s} = \norm{\angles{\xi}^s \widehat f(\xi)}_{L^2_\xi}.
$$
Here $\angles{\cdot} = 1 + \abs{\cdot}$. For $a,\alpha \in \R$, let $X^{a,\alpha}_\pm$ and $H^{a,\alpha}$ be the completions of $\mathcal S(\R^{1+1})$ with respect to
\begin{align*}
  \norm{u}_{X^{a,\alpha}_\pm} &= \bignorm{\angles{\xi}^a
\angles{\tau\pm \xi}^\alpha \widetilde u(\tau,\xi)}_{L^2_{\tau,\xi}},
  \\
  \norm{u}_{H^{a,\alpha}} &= \bignorm{\angles{\xi}^a \angles{\abs{\tau} -
\abs{\xi}}^\alpha \widetilde u(\tau,\xi)}_{L^2_{\tau,\xi}}.
   \end{align*}
See \cite{Selberg:2006b} for more details about these spaces. Finally, if $X,Y,Z$ are normed function spaces, we use the notation
$$
  X \cdot Y \hookrightarrow Z
$$
to mean that $\norm{uv}_Z \lesssim \norm{u}_X \norm{v}_Y$.

\section{Preliminaries}

The Dirac operator $\alpha D_x$ has Fourier symbol $\alpha \xi$, whose eigenvalues are $\pm \xi$. The eigenspace projections are
$$
  P_{\pm} = \frac12
  \begin{pmatrix}
    1 & \pm1  \\
    \pm1 & 1
  \end{pmatrix}.
$$
Following \cite{Selberg:2006b}, Pecher used instead the ordering $\pm \abs{\xi}$ of the eigenvalues, yielding nonconstant projections (with our choice of $\alpha,\beta$)
$$
  \pi_\pm(\xi) = \frac12
  \begin{pmatrix}
    1 & \pm \sgn \xi  \\
    \pm \sgn \xi & 1
  \end{pmatrix}.
$$
The fact that our projections are constant simplifies the argument considerably.

We now write $\psi = \psi_+ + \psi_-$, where
$$
  \psi_+ = P_+ \psi = \frac12 \begin{pmatrix}
     \psi_1 + \psi_2 \\
     \psi_1 + \psi_2 
  \end{pmatrix},
  \quad
  \psi_- = P_- \psi = \frac12 \begin{pmatrix}
     \psi_1 - \psi_2 \\
     \psi_2 - \psi_1 
  \end{pmatrix}.
$$
Applying $P_{\pm}$ on both sides of the first equation in
(\ref{DKG1}), and using the identities $\alpha = P_+ -  P_-$, $P_{\pm}^2 = P_{\pm}$ and $P_{\pm}P_{\mp}=0$, (\ref{DKG1}) is rewritten as
\begin{equation}\label{DKG2}
\left\{
\begin{aligned}
  &(D_t +  D_x ) \psi_+ = P_+(\phi \beta\psi),
  \\
  &(D_t -  D_x ) \psi_- = P_-(\phi \beta\psi),
  \\
  &\square   \phi = -\innerprod{\beta\psi
}{\psi}_{\C^2}.
\end{aligned}
\right.
\end{equation}
We iterate in the spaces
$$
  \psi_+ \in X_+^{s,\sigma},
  \quad
  \psi_- \in X_-^{s,\sigma},
  \quad
  (\phi,\partial_t \phi) \in H^{r,\rho} \times H^{r-1,\rho},
$$
where
$$
  \frac{1}{2} < \sigma, \rho \le 1
$$
will be chosen depending on $r,s$. By a standard argument (see \cite{Selberg:2006b} for details) Theorem \ref{Thm1} then reduces to
\begin{align}
  \label{BilinearAA}
  \norm{P_\pm(\phi\beta P_{[\pm]}
  \psi)}_{X_{\pm}^{s,\sigma-1+\varepsilon}}
  &\lesssim \norm{\phi}_{H^{r,\rho}}
  \norm{\psi}_{X_{[\pm]}^{s,\sigma}},
  \\
  \label{BilinearBB}
  \norm{\innerprod{\beta P_{[\pm]}\psi}{P_{\pm}\psi'}_{\C^2}}_{H^{r-1,\rho-1+\varepsilon}}
  &\lesssim
  \norm{\psi}_{{X_{[\pm]}^{s,\sigma}}}
  \norm{\psi'}_{X_{\pm}^{s,\sigma}},
\end{align}
where $\pm$ and $[\pm]$ denote independent signs and $\varepsilon > 0$ is sufficiently small; the introduction of the parameter $\varepsilon$ is a technical detail needed in the time localized linear estimates (see \cite[Lemmas 5 and 6]{Selberg:2006b}).

But by a duality argument introduced in \cite{Selberg:2006b}, estimate \eqref{BilinearAA} is in fact equivalent to
\begin{equation}\label{BilinearAAA}\tag{$\text{\ref{BilinearAA}}'$}
\norm{\innerprod{\beta P_{[\pm]}\psi}{P_{\pm}\psi'}_{\C^2}}_{H^{-r,-\rho}}
  \lesssim
\norm{\psi}_{{X_{[\pm]}^{s,\sigma}}}
  \norm{\psi'}_{X_{\pm}^{-s,1-\sigma-\varepsilon}}.
\end{equation}
The advantage of this formulation is that, like \eqref{BilinearBB}, it contains the bilinear form $\innerprod{\beta P_{[\pm]}\psi}{P_{\pm}\psi'}$, which turns out be a null form:
With our choice of projections, this comes out very easily, since by the self-adjointness, idempotency and orthogonality of the $P_\pm$, as well as the identity $P_\pm \beta = \beta P_\mp$, we see that
$$
  \innerprod{\beta P_+\psi}{P_+\psi'}_{\C^2}
  = \innerprod{\beta P_-\psi}{P_-\psi'}_{\C^2}
  = 0.
$$
As a result, ($\text{\ref{BilinearAA}}'$) and \eqref{BilinearBB} can be reduced to
\begin{align}
  \label{BilinearAAAA}
  \norm{u\bar{v}}_{H^{-r,-\rho}}
  &\lesssim
  \norm{u}_{{X_{+}^{s,\sigma}}}
  \norm{v}_{X_{-}^{-s,1-\sigma-\varepsilon}}.
  \\
  \label{BilinearBBB}
  \norm{u \bar{v}}_{H^{r-1,\rho-1+\varepsilon}}
  &\lesssim
  \norm{u}_{{X_{+}^{s,\sigma}}}
  \norm{v}_{X_{-}^{s,\sigma}},
\end{align}
where $u,v$ are $\C$-valued and $\bar v$ denotes the complex conjugate.
The crucial point to note here is the difference in signs on the right, due to the null structure; if we had two equal signs, then the estimates would fail at the regularity prescribed in Theorem \ref{Thm1} (cf.\ the conditions in Theorem \ref{Thm2} below). There are two key reasons why things are better when the signs are different: The first reason is the algebraic constraint given in Lemma \ref{Lemma3} below, which is the analogue, in the current setting, of Lemma 7 in \cite{Selberg:2006b}; the second reason is the bilinear estimate given in Lemma \ref{WaveLemma} below.

\begin{lemma}\label{Lemma3}
Define, for $\tau,\lambda,\xi,\eta \in \R$,
$$
 \Gamma=\abs{\tau}-\abs{\xi}, \quad
 \Theta_+=\lambda+\eta, \quad
 \Sigma_-=\lambda-\tau-(\eta-\xi).
$$
Then
$$
 \min( \abs{\eta}, \abs{\eta-\xi})
 \le \frac32 \max\left( \abs{\Gamma}, \abs{\Theta_+}, \abs{\Sigma_-} \right).
$$
\end{lemma}

\begin{proof} We have
$$
  \Gamma =
  \begin{cases}
     \Theta_+ - \Sigma_- - (2\eta-\xi+\abs{\xi}) \quad &\text{if $\tau \ge 0$},
     \\
     -\Theta_+ + \Sigma_- + (2\eta-\xi-\abs{\xi}) \quad &\text{if $\tau \le 0$}.
  \end{cases}
$$
and the terms in parentheses equal $2\eta$ or $2(\eta-\xi)$, depending on the sign of $\xi$. Therefore, $2\min(\abs{\eta},\abs{\eta-\xi}) \le \abs{\Gamma} + \abs{\Theta_+} + \abs{\Sigma_-}$.
\end{proof}

This lemma is applied in tandem with the following product law for the Wave-Sobolev spaces $H^{a,\alpha}$. The sufficiency of the conditions \eqref{abc1} and \eqref{abc2} in the following theorem can easily be deduced from \cite[Proposition A.1]{Selberg:2002b}, but here we also prove necessity, up to endpoints; see Section \ref{Thm2Proof}.

\begin{theorem}\label{Thm2}
Suppose $a,b,c \in \R$, \,$\alpha, \beta, \gamma \ge 0$
and $\alpha + \beta + \gamma > \frac{1}{2}$. Then
\begin{equation}\label{Embed}
  H^{a, \alpha} \cdot H^{b, \beta} \embeds H^{-c,
-\gamma},
\end{equation}
provided that
\begin{gather}
  \label{abc1}
  a + b + c > \frac12,
  \\
  \label{abc2}
  a + b \ge 0, \quad a + c \ge 0, \quad b + c \ge 0.
\end{gather}
Furthermore, these conditions are sharp up to
equality, in the sense that if \eqref{Embed} holds, then \eqref{abc2} must hold, and \eqref{abc1} must hold with $\ge$.
\end{theorem}

\begin{remark} The above product law is analogous to the one for the standard Sobolev spaces, which in 1d reads $\norm{fg}_{H^{-c}} \lesssim \norm{f}_{H^a} \norm{g}_{H^b}$, with the same conditions on $a,b,c$ as in the above theorem.
\end{remark}

The algebraic constraint (Lemma \ref{Lemma3}) and the product law for Wave-Sobolev spaces are enough to prove the result of Pecher, but to improve on that result, we use also the following bilinear space-time estimate for 1d free waves, where again the different signs are crucial.

\begin{lemma}\label{WaveLemma}
Suppose $u,v$ solve
$$
\begin{alignedat}{2}
  &(D_t + D_x ) u = 0,& \qquad &u(0,x) = f(x),
  \\
  &(D_t - D_x ) v = 0,& \qquad &v(0,x) = g(x),
\end{alignedat}
$$
where $f,g \in L^2(\R)$. Then
$$
  \norm{uv}_{L^2(\R^{1+1})} \le \sqrt{2} \norm{f}_{L^2} \norm{g}_{L^2}.
$$
\end{lemma}

\begin{proof}
We have $\widetilde u(\tau,\xi) = \delta(\tau+\xi) \widehat f(\xi)$ and $\widetilde v(\tau,\xi) = \delta(\tau-\xi) \widehat g(\xi)$, so
\begin{align*}
  \widetilde{uv}(\tau,\xi)
  &= \int_{\R^{1+1}} \widetilde u(\lambda,\eta) \widetilde v(\tau-\lambda,\xi-\eta) \, d\lambda \, d\eta
  \\
  &= \int \delta(\tau+2\eta-\xi) \widehat f(\eta) \widehat g(\xi-\eta) \, d\eta
  \\
  &= \widehat f \left( \frac{\xi-\tau}{2} \right) \widehat g \left( \frac{\xi+\tau}{2} \right).
\end{align*}
The claimed estimate now follows from Plancherel's theorem and an obvious change of variables.
\end{proof}

By the transfer principle (see \cite[Lemma 4]{Selberg:2006b}), Lemma \ref{WaveLemma} implies:

\begin{corollary}\label{WaveCorollary} For any $\alpha > 1/2$,
$$
  X_+^{0,\alpha} \cdot X_-^{0,\alpha} \hookrightarrow L^2.
$$
\end{corollary}

Again, this would fail if we had equal signs in the left hand side.

We now have all the tools needed to finish the proof of the main estimates.

\section{Proof of Theorem \ref{Thm1}}\label{Thm1Proof}

\subsection{Proof of \eqref{BilinearAAAA}}

With notation as in Lemma \ref{Lemma3}, the estimate is equivalent to, using Plancherel's theorem,
$$
 \norm{ \int_{\R^2} \frac{ F(\lambda,\eta)
 G(\lambda-\tau,\eta-\xi)d\lambda \, d\eta}
 {\angles{\xi}^{r}\angles{\eta}^{s} \angles{\eta-\xi}^{-s}
  \angles{\Gamma}^{\rho}
  \angles{\Theta_+}^{\sigma}
  \angles{\Sigma_-}^{1-\sigma-\varepsilon}}}_{L^2_{\tau,\xi}}
  \lesssim \norm{F}_{L^2} \norm{G}_{L^2},
$$
for arbitrary $F, G \in L^2(\R^2)$. In view of Lemma \ref{Lemma3} we can add either $\rho$, $\sigma$ or $1-\sigma-\varepsilon$ to the exponent of either the $\angles{\eta}$ weight or the $\angles{\eta-\xi}$ weight, at the expense of giving up one of the ``hyperbolic'' weights $\angles{\Gamma}$, $\angles{\Theta_+}$ or $\angles{\Sigma_-}$. Then we apply Theorem \ref{Thm2}. In fact, since (recall $\rho,\sigma > 1/2$)
$$
  \min(\rho,\sigma,1-\sigma-\varepsilon) = 1-\sigma-\varepsilon,
$$
we can reduce to Theorem \ref{Thm2} with $a,b,c$ as in the first two rows of Table \ref{Table1}. The conditions on $a,b,c$ in Theorem \ref{Thm2} impose the following restrictions:
\begin{gather}
  \label{r1}
  r > \sigma-\frac12+\varepsilon,
  \\
  \label{r2}
  r \ge \abs{s},
  \\
  \label{sigma1}
  \sigma \le 1-\varepsilon.
\end{gather}
Finally, we mention that the hypotheses on $(\alpha,\beta,\gamma)$ in Theorem \ref{Thm2} are indeed satisfied in this situation, as follows from \eqref{sigma1} and the fact that we require
\begin{equation}\label{rho_sigma}
  \frac12 < \rho, \sigma \le 1.
\end{equation}
So we conclude that \eqref{BilinearAAAA} holds provided \eqref{r1}--\eqref{rho_sigma} are verified.

\subsection{Proof of \eqref{BilinearBBB}}

This reduces to
$$
  I := \norm{ \int_{\R^2} \frac{ F(\lambda,\eta)
  G(\lambda-\tau,\eta-\xi)d\lambda \, d\eta}
  {\angles{\xi}^{1-r}\angles{\eta}^{s} \angles{\eta-\xi}^{s}
  \angles{\Gamma}^{1-\rho-\varepsilon}
  \angles{\Theta_+}^{\sigma}
  \angles{\Sigma_-}^{\sigma}}}_{L^2_{\tau,\xi}}
  \lesssim \norm{F}_{L^2} \norm{G}_{L^2}.
$$
We consider two cases, with notation as in Lemma \ref{Lemma3}:

\subsubsection{Case 1: $\max(\abs{\Gamma},\abs{\Theta_+},\abs{\Sigma_-}) \sim \abs{\Gamma}$}

By symmetry we may assume $\abs{\eta} \le \abs{\eta-\xi}$ in $I$. Then either $\abs{\eta} \sim \abs{\eta-\xi}$, or $\abs{\eta} \ll \abs{\eta-\xi} \sim \abs{\xi}$, hence, using Lemma \ref{Lemma3},
$$
  I \lesssim I_1 + I_2
$$
where
\begin{align*}
  I_1 &= \norm{ \int_{\abs{\eta}\sim\abs{\eta-\xi}} \frac{ F(\lambda,\eta)
  G(\lambda-\tau,\eta-\xi)d\lambda \, d\eta}
  {\angles{\xi}^{1-r}\angles{\eta}^{2s+1-\rho-\varepsilon} \angles{\Theta_+}^{\sigma}
  \angles{\Sigma_-}^{\sigma}}}_{L^2_{\tau,\xi}},
  \\
  I_2 &= \norm{ \int \frac{ F(\lambda,\eta)
  G(\lambda-\tau,\eta-\xi)d\lambda \, d\eta}
  {\angles{\eta}^{s+1-\rho-\varepsilon} \angles{\eta-\xi}^{s+1-r} \angles{\Theta_+}^{\sigma}
  \angles{\Sigma_-}^{\sigma}}}_{L^2_{\tau,\xi}}.
\end{align*}
Moreover, if $r > 1$, then we can use $\angles{\xi}^{r-1} \lesssim \angles{\eta}^{r-1} + \angles{\eta-\xi}^{r-1}$ to further reduce $I_1$ to
$$
  I_{1,r>1} = \norm{ \int_{\abs{\eta}\sim\abs{\eta-\xi}} \frac{ F(\lambda,\eta)
  G(\lambda-\tau,\eta-\xi)d\lambda \, d\eta}
  {\angles{\eta}^{2s+1-\rho-\varepsilon+1-r} \angles{\Theta_+}^{\sigma}
  \angles{\Sigma_-}^{\sigma}}}_{L^2_{\tau,\xi}}.
$$
Applying Corollary \ref{WaveCorollary}, we then see that $I_i \lesssim \norm{F}_{L^2} \norm{G}_{L^2}$ ($i=1,2$), provided
\begin{gather}
  \label{r6}
  r \le 1+s,
  \\
  \label{s2}
  s \ge -\frac12+\frac{\rho+\varepsilon}{2}
  \\
  \label{s3}
  s \ge -1+\rho+\varepsilon
  \\
  \label{r7}
  r \le 1+2s+1-\rho-\varepsilon.
\end{gather}
\begin{remark}
If we had applied Theorem \ref{Thm2} here, the last condition would have been replaced by (due to the requirement $a+b+c>1/2$ in Theorem \ref{Thm2})
$$
  r < \frac12+2s+1-\rho-\varepsilon,
$$
which is still sufficient to obtain the result of Pecher. So it is exactly at this point that we gain something more.
\end{remark}

\subsubsection{Case 2: $\max(\abs{\Gamma},\abs{\Theta_+},\abs{\Sigma_-}) \sim \abs{\Theta_+}$ or $\abs{\Sigma_-}$}

Then by Lemma \ref{Lemma3} we reduce to Theorem \ref{Thm2} with $a,b,c$ as in the last row of Table \ref{Table1}, and $(\alpha,\beta,\gamma)=(0,\sigma,1-\rho-\varepsilon)$ or $(\sigma,0,1-\rho-\varepsilon)$. The conditions on $a,b,c,\alpha,\beta,\gamma$ in Theorem \ref{Thm2} yield the restrictions
\begin{gather}
  \label{r3}
  r < \frac12 + \sigma + 2s,
  \\
  \label{r4}
  r \le 1+s,
  \\
  \label{s1}
  s \ge -\frac{\sigma}{2}
  \\
  \label{rho1}
  \rho \le 1-\varepsilon.
\end{gather}
Note that \eqref{r4} is the same as \eqref{r6}.

\begin{table}
\caption{Exponents used in Theorem \ref{Thm2}}
\label{Table1}
\def\arraystretch{1.3}
\begin{center}
\begin{tabular}{|c|c|c|}
  \hline
  $a$ & $b$ & $c$
  \\
  \hline
  \hline
  $s$ & $-s+1-\sigma-\varepsilon$ & $r$ 
  \\
  \hline
  $s+1-\sigma-\varepsilon$ & $-s$ & $r$
  \\
  \hline
  $s+\sigma$ & $s$ & $1-r$ 
  \\
  \hline
\end{tabular}
\end{center}
\end{table}

\medskip\noindent
We conclude that \eqref{BilinearBBB} holds if \eqref{r6}--\eqref{rho1} are satisfied.

\subsection{Conlusion of the proof}

It only remains, given $(s,r)$ satisfying the hypotheses
\begin{equation}\label{hypotheses}
  s > -\frac{1}{4}, \qquad r>0,
  \qquad
  \abs{s}
  \le r \le 1+s
\end{equation}
of Theorem \ref{Thm1}, to choose $\rho,\sigma,\varepsilon$ in such a way that the constraints \eqref{r1}--\eqref{rho1} are all satisfied. We shall need the fact that \eqref{hypotheses} implies
\begin{equation}\label{hypotheses'}
  r < 3/2 + 2s.
\end{equation}

Clearly, we get the best results by choosing $\rho$ and $\varepsilon$ as small as possible, so let us set
$$
  \rho = \frac12+\varepsilon,
$$
where $\varepsilon > 0$ will be chosen sufficiently small. Note that \eqref{rho1} is satisfied provided $\varepsilon \le 1/4$. Condition \eqref{s2} becomes
$$
  s \ge -\frac14 + \varepsilon,
$$
which is compatible with the assumption $s > -1/4$ in Theorem \ref{Thm1}; \eqref{s3} and \eqref{s1} are weaker than \eqref{s2}, so they are also satisfied. Condition \eqref{r7} becomes
$$
  r \le \frac32 + 2s - 2\varepsilon,
$$
in accordance with \eqref{hypotheses'}.

The only remaining conditions are \eqref{r1}, \eqref{sigma1} and \eqref{r3} (as well as \eqref{rho_sigma}, which requires $\sigma > 1/2$), and these conditions can be summed up as follows:
\begin{gather*}
  \frac12 < \sigma \le 1-\varepsilon,
  \\
  \sigma-\frac12+\varepsilon < r < \sigma + \frac12 + 2s.
\end{gather*}
Since $0 < r < 3/2+2s$, by \eqref{hypotheses} and \eqref{hypotheses'}, it is clear that we can find $\sigma$ and $\varepsilon > 0$ such that the last two conditions are satisfied. This completes the proof of Theorem \ref{Thm1}.

\section{Counterexamples}\label{Optimality}

Here we prove the optimality, except for the endpoint $(s,r) = (0,0)$, of the conditions on $s$ and $r$ in Theorem \ref{Thm1}, as far as iteration in the Bourgain-Klainerman-Machedon spaces $X_\pm^{s,\sigma}$, $H^{r,\rho}$ is concernced. To be precise, we prove:

\begin{theorem}\label{Thm3} \ 
\begin{enumerate}
\item
The estimate \eqref{BilinearBB} fails (for every choice of $1/2 < \sigma, \rho \le 1$ and $\varepsilon > 0$) if $s \le - 1/4$ or $r > 1+s$.
\item
The estimate \eqref{BilinearAAA}, hence also \eqref{BilinearAA}, fails (for every choice of $1/2 < \sigma, \rho \le 1$ and $\varepsilon > 0$) if $r < \abs{s}$.
\end{enumerate}
\end{theorem}

More generally, we prove:
\begin{theorem}\label{Thm4}
Let $a,b,c,\alpha,\beta,\gamma \in \R$. If the 2-spinor estimate 
$$
  \norm{\innerprod{\beta P_+\psi}{P_-\psi'}_{\C^2}}_{H^{-c,-\gamma}}
  \lesssim
  \norm{\psi}_{{X_+^{a,\alpha}}}
  \norm{\psi'}_{X_-^{b,\beta}},
$$
holds, then:
\begin{gather}
  \label{Cond1}
  a+b+\min(\alpha,\beta,\gamma) \ge 0,
  \\
  \label{Cond2}
  a+b+c+\min(\alpha,\beta) \ge \frac12
  \\
  \label{Cond4}
  a+b+c+\gamma \ge 0,
  \\
  \label{Cond3}
  \min(a,b)+c \ge 0.
\end{gather}
\end{theorem}

\subsection{Proof of Theorem \ref{Thm3}}
We apply Theorem \ref{Thm4}. For part (a) we take $(a,b,c) = (s,s,1-r)$ and $(\alpha,\beta,\gamma) = (\sigma,\sigma,1-\rho-\varepsilon)$. Then \eqref{Cond1} gives
the necessary condition $2s + 1-\rho-\varepsilon \ge 0$, i.e., $s \ge -1/2+(\rho+\varepsilon)/2 > -1/4$, where the last inequality holds since $\rho > 1/2$. Moreover, \eqref{Cond3} gives the necessary condition $s+1-r \ge 0$. This proves part (a).

To prove part (b), take $(a,b,c) = (s,-s,r)$ and $(\alpha,\beta,\gamma) = (\sigma,1-\sigma-\varepsilon,\rho)$. Then \eqref{Cond3} implies $-\abs{s} + r \ge 0$.

Note that we only used \eqref{Cond1} and \eqref{Cond3} to prove Theorem \ref{Thm3}.

\subsection{Proof of Theorem \ref{Thm4}}

The following counterexamples are adapted from those for the 2d case in \cite{Selberg:2006c}, and depend on a large, positive parameter $L$ going to infinity. We choose intervals $A, B, C \subset
\R$, depending on $L$, with the property
\begin{equation}\label{ABCproperty}
  \eta \in A, \,\, \xi \in C \implies \eta-\xi \in B.
\end{equation}
We shall denote by $\abs{A}$ the length of the interval $A$.

We shall set
\begin{equation}\label{bad_spinors}
  \psi(t,x) = u(t,x) {1\choose1}, \quad
  \psi'(t,x) = v(t,x) {1\choose-1},
\end{equation}
where $u,v : \R^{1+1} \to \C$ are defined on the Fourier transform side by
\begin{equation}\label{u_v}
  \widetilde u(\lambda,\eta) = \mathbf{1}_{\lambda+\eta= O(1)}
  \mathbf{1}_{\eta \in A},
  \quad
  \widetilde v(\lambda-\tau,\eta-\xi) = \mathbf{1}_{\lambda-\tau+\eta-\xi = O(1)}
  \mathbf{1}_{\eta -\xi \in B},
\end{equation}
and $A$, $B$ remain to be chosen. Here $\mathbf{1}_{(\cdot)}$ stands for the indicator function of the set determined by the condition in the subscript. Then
\begin{equation}\label{innerprod}
  \innerprod{\beta P_+\psi}{P_-\psi'}_{\C^2} = \innerprod{\beta \psi}{\psi'}_{\C^2} = 2 u\bar v,
\end{equation}
so in fact it suffices to find counterexamples to
\begin{equation}\label{bad_uv}
  \norm{u\bar v}_{H^{-c,-\gamma}}
  \lesssim
  \norm{u}_{{X_+^{a,\alpha}}}
  \norm{v}_{X_-^{b,\beta}}.
\end{equation}
Each counterexample will be of the form
\begin{equation}\label{Product1}
  \frac{\norm{u\bar{v}}_{H^{-c,-\gamma}}}{\norm{u}_{X_+^{a,\alpha}}\norm{v}_{X_-^{b,\beta}}} \gtrsim \frac{1}{L^{\delta(a,b,c,\alpha,\beta,\gamma)}},
\end{equation}
which leads to the necessary condition $\delta(a,b,c,\alpha,\beta,\gamma) \ge 0$.

Observe that
\begin{align*}
  &\norm{u\bar{v}}_{H^{-c,-\gamma}}
  \\
  &= \norm{\int_{\R^{1+1}}
  \frac{1}{\angles{\xi}^{c}\angles{\abs{\tau}-\abs{\xi}}^{\gamma}}
  \mathbf{1}_{\left\{
  \scriptstyle\eta \in A
  \atop
  \scriptstyle\lambda+\eta= O(1)
  \right\}
  }
  \mathbf{1}_{\left\{
    \scriptstyle\eta-\xi \in B
    \atop
  \scriptstyle\lambda-\tau+\eta-\xi = O(1)
  \right\}
  }
  \, d\lambda \, d\eta}_{L^2_{\tau,\xi}}
  \\
  &\ge I := \norm{\int_{\R^{1+1}}
  \frac{1}{\angles{\xi}^{c}\angles{\abs{\tau}-\abs{\xi}}^{\gamma}}
  \mathbf{1}_{\left\{
  \scriptstyle\eta \in A
  \atop
  \scriptstyle\lambda+\eta = O(1)
  \right\}
  }
  \mathbf{1}_{\left\{
  \scriptstyle\xi \in C
  \atop
  \scriptstyle\tau+\xi = O(1)
  \right\}
  }
  \, d\lambda \, d\eta}_{L^2_{\tau,\xi}},
\end{align*}
where to get the last inequality we restrict the $L^2$ norm to $\tau+\xi=O(1)$, $\xi \in C$, make use of \eqref{ABCproperty}, and note that
$$
  \lambda+\eta = O(1), \,\, \tau+\xi = O(1)
  \implies \lambda-\tau+\eta-\xi = O(1).
$$
Note also that $\bigabs{ \abs{\tau}-\abs{\xi}} \le \abs{\tau + \xi} = O(1)$. (So to get counterexamples involving $\gamma$, we shall later have to modify $I$).

\subsubsection{Necessity of \eqref{Cond1} when $\min(\alpha,\beta,\gamma) = \alpha$ or $\beta$}
Define
$$
  A = [L-1/2,L+1/2],
  \quad
  B = [L-1,L+1],
  \quad
  C = [-1/2,1/2].
$$
Then $\abs{\xi} = O(1)$, $\fixedabs{\eta} \sim L$, $\fixedabs{\eta-\xi} \sim L$ and
$$
  \lambda-\tau-(\eta-\xi) = \lambda-\tau+(\eta-\xi) - 2(\eta-\xi) \sim L,
$$
hence
$$
  I \sim \abs{A} \abs{C}^{1/2},
  \quad
  \norm{u}_{X_+^{a,\alpha}} \sim L^{a} \abs{A}^{1/2},
  \quad
  \norm{v}_{X_-^{b,\beta}} \sim L^{b+\beta} \abs{B}^{1/2}.
$$
But $\abs{A}, \abs{B}, \abs{C} \sim 1$, so \eqref{Product1} holds with $\delta(a,b,c,\alpha,\beta,\gamma) = a+b+\beta$, which gives the necessary condition $a+b+\beta \ge 0$. By symmetry, we must also have $a+b+\alpha \ge 0$.

\subsubsection{Necessity of \eqref{Cond2}}\label{Counterex2}
Set
$$
  A = [L/4,L/2],
  \quad
  B = [L/2,3L/2],
  \quad
  C = [-L,-L/2].
$$
Then $\fixedabs{\eta}, \abs{\xi}, \fixedabs{\eta-\xi} \sim L$ and (as above)
$\lambda-\tau-(\eta-\xi) \sim L$, so
$$
  I \sim \frac{\abs{A}\abs{C}^{1/2}}{L^{c}},
  \quad
  \norm{u}_{X_+^{a,\alpha}} \sim L^{a} \abs{A}^{1/2},
  \quad \norm{v}_{X_-^{b,\beta}}
  \sim L^{b+\beta} \abs{B}^{1/2}.
$$
Since $\abs{A}, \abs{B}, \abs{C} \sim L$, we conclude that
\eqref{Product1} holds with $\delta(a,b,c,\alpha,\beta,\gamma) = a+b+c+\beta-1/2$, proving the necessity of $a+b+c+\beta \ge 1/2$. By symmetry, we also need $a+b+c+\alpha \ge 1/2$

\subsubsection{Necessity of \eqref{Cond3}}\label{Counterex3} Here we set
$$
  A = C = [L-1/2,L+1/2],
  \quad
  B = [-1,1].
$$
Then $\abs{\xi} \sim L$, $\fixedabs{\eta} \sim L$, $\fixedabs{\eta-\xi} = O(1)$ and
$$
  \lambda-\tau-(\eta-\xi) =\lambda-\tau+(\eta-\xi)-2(\eta-\xi) = O(1),
$$
so
$$
  I \sim \frac{\abs{A}\abs{C}^{1/2}}{L^{c}},
  \qquad
  \norm{u}_{X_+^{a,\alpha}} \sim L^a \abs{A}^{1/2},
  \qquad \norm{v}_{X_-^{b,\beta}}
  \sim \abs{B}^{1/2}.
$$
But $\abs{A}, \abs{B}, \abs{C} \sim 1$, hence \eqref{Product1} holds with $\delta(a,b,c,\alpha,\beta,\gamma) = a+c$, proving necessity of $a+c \ge 0$. By symmetry, $b+c \ge 0$ is also necessary.

\subsubsection{Necessity of \eqref{Cond1} when $\min(\alpha,\beta,\gamma) = \gamma$}\label{GammaCase}
Set
$$
  A = [L-1,L+1],
  \quad
  B = [L-2,L+2],
  \quad
  C = [-1,1].
$$
Again we use \eqref{bad_spinors}, with $u$ as in \eqref{u_v}, but we change $v$ to:
$$
  \widetilde v(\lambda-\tau,\eta-\xi) = \mathbf{1}_{\lambda-\tau-(\eta-\xi) = O(1)}
  \mathbf{1}_{\eta -\xi \in B}.
$$
Since \eqref{innerprod} is unchanged, it suffices to disprove \eqref{bad_uv}, but now, in view of the modification of $v$,
\begin{align*}
  &\norm{u\bar{v}}_{H^{-c,-\gamma}}
  \\
  &= \norm{\int_{\R^{1+1}}
  \frac{1}{\angles{\xi}^{c}\angles{\abs{\tau}-\abs{\xi}}^{\gamma}}
  \mathbf{1}_{\left\{
  \scriptstyle\eta \in A
  \atop
  \scriptstyle\lambda+\eta= O(1)
  \right\}
  }
  \mathbf{1}_{\left\{
    \scriptstyle\eta-\xi \in B
    \atop
  \scriptstyle\lambda-\tau-(\eta-\xi) = O(1)
  \right\}
  }
  \, d\lambda \, d\eta}_{L^2_{\tau,\xi}}
  \\
  &\ge I := \norm{\int_{\R^{1+1}}
  \frac{1}{\angles{\xi}^{c}\angles{\abs{\tau}-\abs{\xi}}^{\gamma}}
  \mathbf{1}_{\left\{
  \scriptstyle\eta \in A
  \atop
  \scriptstyle\lambda+\eta = O(1)
  \right\}
  }
  \mathbf{1}_{\left\{
  \scriptstyle\xi \in C
  \atop
  \scriptstyle\tau+2L = O(1)
  \right\}
  }
  \, d\lambda \, d\eta}_{L^2_{\tau,\xi}},
\end{align*}
where in the last step we restrict the $L^2$ norm to the region
$\tau+2L = O(1)$, $\xi \in C$, make use of \eqref{ABCproperty}, and note that
$$
  \lambda-\tau-(\eta-\xi) = (\lambda+\eta) + 2(L-\eta) - (\tau+2L) + \xi = O(1),
$$
since each term is $O(1)$. So now $\abs{\xi} = O(1)$, $\abs{\eta} \sim L$, $\abs{\eta-\xi} \sim L$, and $\bigabs{\abs{\tau}-\abs{\xi}} \sim L$, hence
$$
  I \sim \frac{\abs{A}\abs{C}^{1/2}}{L^{\gamma}},
  \quad
  \norm{u}_{X_+^{a,\alpha}} \sim L^{a} \abs{A}^{1/2},
  \quad \norm{v}_{X_-^{b,\beta}}
  \sim L^{b} \abs{B}^{1/2}.
$$
Since $\abs{A}, \abs{B}, \abs{C} \sim 1$, \eqref{Product1} holds with $\delta(a,b,c,\alpha,\beta,\gamma) = a+b+\gamma$.

\subsubsection{Necessity of \eqref{Cond4}}
Here we use the same $u,v$ as in subsection \ref{GammaCase}. Set
$$
  A = [L-1,L+1],
  \quad
  B = [2L-2,2L+2],
  \quad
  C = [-L-1,-L+1].
$$
Then as in subsection \ref{GammaCase}, we have $\norm{u\bar{v}}_{H^{-c,-\gamma}} \ge I$, with the only difference that the condition $\tau+2L=O(1)$ in $I$ has been replaced by $\tau+3L=O(1)$, for then we can write
$$
  \lambda-\tau-(\eta-\xi) = (\lambda+\eta) + 2(L-\eta) - (\tau+3L) + (\xi+L) = O(1),
$$
each term being $O(1)$. So
$\abs{\xi}, \abs{\eta}, \abs{\eta-\xi} \sim L$, and $\bigabs{\abs{\tau}-\abs{\xi}} \sim L$, hence
$$
  I \sim \frac{\abs{A}\abs{C}^{1/2}}{L^{c+\gamma}},
  \quad
  \norm{u}_{X_+^{a,\alpha}} \sim L^{a} \abs{A}^{1/2},
  \quad \norm{v}_{X_-^{b,\beta}}
  \sim L^{b} \abs{B}^{1/2}.
$$
Since $\abs{A}, \abs{B}, \abs{C} \sim 1$, \eqref{Product1} holds with $\delta(a,b,c,\alpha,\beta,\gamma) = a+b+c+\gamma$.

\section{Proof of Theorem \ref{Thm2}}\label{Thm2Proof}

In this section we fix $\alpha, \beta, \gamma \ge 0$ satisfying $\alpha + \beta + \gamma > 1/2$. We shall say that a triple $(a,b,c)$ of real numbers is \emph{admissible} if the embedding \eqref{Embed} holds, i.e., if the bilinear estimate
\begin{equation}\label{Product2}
  \norm{uv}_{ H^{-c,- \gamma} }
  \lesssim \norm{u}_{ H^{a,\alpha} } \norm {v}_{ H^{b,\beta} }
\end{equation}
holds.

First, assume that conditions \eqref{abc1} and \eqref{abc2} are satisfied. If $a, b, c \ge 0$, then $(a,b,c)$ is admissible, as proved in \cite[Proposition A.1]{Selberg:2002b}. It remains to consider the case where $(a,b,c)$ contains a negative number. But in view of \eqref{abc2}, at most one of the numbers $a,b,c$ can be negative, and by symmetry it suffices to consider the case $a<0$, say. In that case we can write $\angles{\xi}^{-a} \lesssim \angles{\eta}^{-a} + \angles{\eta+ \xi}^{-a}$, thus reducing to the triples $(0,a+b,c)$ or $(0,b,a+c)$, which contain no negative numbers, hence are admissible, as noted above.

It remains to prove necessity of \eqref{abc1} (up to equality) and \eqref{abc2}. But in fact, the counterexample constructed in subsection \ref{Counterex2} gives, with $u,v$ as in \eqref{u_v}
\begin{equation}\label{counterex}
  \frac{\norm{u\bar{v}}_{H^{-c,-\gamma}}}{\norm{u}_{H^{a,\alpha}}\norm{v}_{H^{b,\beta}}} \gtrsim \frac{1}{L^{\delta(a,b,c,\alpha,\beta,\gamma)}},
\end{equation}
with $\delta(a,b,c,\alpha,\beta,\gamma) = a + b + c - 1/2$, proving necessity of $a+b+c \ge 1/2$, which is \eqref{abc1} with $\ge$. The fact that we have a conjugate here, but not in \eqref{Product2}, is irrelevant, since $\norm{v}_{H^{b,\beta}} = \norm{\bar v}_{H^{b,\beta}}$.

Similarly, the counterexample from \ref{Counterex3} gives \eqref{counterex} with
$\delta(a,b,c,\alpha,\beta,\gamma) = a + c$, proving necessity of $a+c \ge 0$. By duality and symmetry in \eqref{Product2}, we then get also $a+b \ge 0$ and $b+c \ge 0$, so we have proved necessity of \eqref{abc2}.

\bibliographystyle{amsplain}
\bibliography{/Users/sselberg/Mathematics/Bibliography/mybibliography}

\end{document}